 \documentclass[final,3p,times]{elsarticle}

\usepackage{setspace}

\usepackage{hyperref}
 \usepackage{graphicx}
\graphicspath{{figs/}}
 \usepackage{color}
\usepackage{epstopdf}
\usepackage{tabularx}
\usepackage{subfig}
\usepackage[font=small,labelfont=bf]{caption}

\usepackage{amssymb}
\usepackage{amsmath}
\usepackage{multirow}
\usepackage{array}

  \usepackage{tabu}
\journal{ }

\begin{document}
 	\doublespacing

\begin{frontmatter}



\title{3D architected isotropic materials with tunable stiffness  and buckling strength }


 \author{Fengwen~Wang\corref{cor1}}
   \ead{fwan@mek.dtu.dk}
     \author{O. Sigmund}

  \address{Department of Mechanical Engineering, Technical University of Denmark, \\ Nils Koppels All{\'e}, Building 404, 2800 Kgs. Lyngby, Denmark}
\cortext[cor1]{Corresponding author. Tel. :+45 4525 4266;  Fax: +45 4593 1475.}
\begin{abstract}
This paper presents a class of  3D single-scale isotropic materials with tunable stiffness and buckling strength obtained via topology optimization and subsequent shape optimization.  Compared to stiffness-optimal closed-cell plate material, the material class reduces the Young's modulus to a range from $79\%$ to $58\%$, but  improves the uniaxial buckling strength to a range from  $180\%$ to $767\%$. Based on small deformation theory, material stiffness is evaluated using the homogenization method. Buckling strength under a given macroscopic stress state is estimated using linear buckling analysis with  Block-Floquet boundary conditions to capture both short   and long wavelength buckling modes.  The 3D isotropic single-scale materials with tunable properties are  designed using topology optimization,  and  are then further simplified using  shape optimization.  Both topology and shape optimized results demonstrate that material buckling strength can be significantly enhanced by hybrids between truss and variable thickness plate structures.

\end{abstract}


\end{frontmatter}

\section{Introduction}
 
Materials with extreme mechanical properties are highly attractive for many applications. Among them, stiffness and strength are fundamental for determining material load-bearing capability.  Stiffness accounts for the ability to resist deformation while strength measures the ultimate load-carrying capability. 
 
Many studies have been devoted to exploring materials with optimal stiffness~ \cite{Hashin1962, Francfort1986,Milton1986,Sigmund1994,Sigmund2000,Berger2017,Tancogne-Dejean2018a,Wang2019}.  It has been shown that isotropic stiffness-optimal plate materials can achieve the Hashin–Shtrikman upper bounds on isotropic elastic stiffness and show up to three times the stiffness of the isotropic stiffness-optimal truss materials in the low volume fraction limit~\cite{Berger2017,Tancogne-Dejean2018a,Christensen1986,Sigmund2016}. The stiffness advantage is attributed to the  multiaxial stiffness offered by the constituent plates while bars in the truss materials only offer axial stiffness.  However, a recent study has shown that isotropic stiffness-optimal truss material is superior from a bucking strength point of view compared to the isotropic plate counterparts for the same volume fraction  when the volume fraction is below 31\%~\cite{Andersen2020} because of higher bending stiffness associated with the bars than the constituent plates.  It is for example shown that the isotropic stiffness-optimal truss material possesses $48\%$ higher uniaxial buckling strength and  $52\%$ lower Young's modulus than the isotropic plate material for a volume fraction of $20\%$.

 Material geometry strongly dictates material properties. Novel materials with exotic properties have been achieved through careful tailoring of material geometries via different design approaches.  Among them, topology optimization methods~\cite{Bendsoe2003} have been proven  powerful tools in designing novel materials ranging from mechanical properties, such as optimal  stiffness~\cite{Sigmund1994,Andreassen2014},  auxetic behavior~\cite{Sigmund1994,Andreassen2014,Wang2018},  to acoustic and optical properties~\cite{Christiansen2016,Wang2018b}.   Regarding elastic stability, material buckling failure may develop at different scales spanning from highly localized short wavelength modes   to long wavelength modes.   Previous numerical studies have employed homogenization methods assuming  separation of scales~\cite{Guedes1990} and Bloch-Floquet theory for detecting short and long wavelength buckling~\cite{Geymonat1993}. A general methodology for characterizing material strength due to bifurcation failure was proposed in  \cite{Triantafyllidis1993}.   Assuming small strains and ignoring material and geometric nonlinearities, topology optimization of material strength was first studied in~ \cite{Neves2002}, where only cell-periodic buckling modes were taken into account. Later, this work was extended to cover both local and global modes via the Bloch-Floquet theory~\cite{Neves2002d}. More recently, 2D materials with enhanced buckling strength~\cite{Thomsen2018} have been systematically designed using topology optimization methods for different macro-level stress situations based on the homogenization theory and linear buckling analysis (LBA) with  Bloch-Floquet boundary conditions.  It was shown that optimized first-order hierarchical materials outperform their non-hierarchical counterparts in terms of buckling strength at the cost of slightly decreased stiffness.  Further material evaluations considering both geometrical and material nonlinearity have proven that the superior buckling strength of the optimized hierarchical materials  compared to reference materials also hold for finite structures~\cite{Wang2020} and  geometrically nonlinear modeling~\cite{Bluhm2020}.  Hence, a buckling optimization approach assuming small strains has been demonstrated to be efficient and practically relevant, even for nonlinear structures.  
 
 This study extends the work in 2D material design with enhanced buckling strength in \citep{Thomsen2018} to 3D  material design with tunable  stiffness and buckling  response utilizing a flexible framework for large scale topology optimization based on the Portable, Extensible Toolkit for Scientific Computation (PETSc) \cite{Balay2016,Aage2015,Wang2018} and the Scalable Library for Eigenvalue Problem Computations (SLEPc)~\cite{Hernandez2005}.     Considered materials are constrained to be cubic symmetric and elastically isotropic.  The homogenization method is employed to characterize the effective material properties, and LBA, together with  Bloch-Floquet boundary conditions evaluated over the boundaries of the irreducible Brillouin zone (IBZ)~\cite{Brillouin1953},  is employed to evaluate material buckling strength. The optimization problem for designing materials with tunable stiffness and buckling strength is formulated to minimize the weighted stiffness and buckling response.   3D materials are designed to achieve tunable stiffness and buckling response by assigning different weight factors for stiffness and buckling strength under uniaxial compression. Moreover, inspired by the topology optimized material configurations, a subsequent feature-based shape optimization approach is employed to simplify material geometries.  In the feature-based approach,  material architectures are parametrized using several hollow and one solid super-ellipsoids~\cite{Wang2018,Wein2020}. 
 
 The paper is organized as follows. Section 2 presents finite element formulations to evaluate stiffness and buckling strength, and formulates the optimization problem for designing materials with tunable properties.  Section 3  first validates the proposed optimization formulation by optimizing a material microstructure for maximum buckling strength under hydrostatic compression.     Topology optimized single-length scale material microstructures with enhanced stiffness and strength under uniaxial compression are then systematically designed. Inspired by the optimized microstructure configurations, a shape optimization scheme is proposed to simplify the optimized microstructure geometries further.    The paper ends with the conclusions in Section 4.

\section{Optimization formulation of 3D  material design with enhanced stiffness and buckling strength}
  This section presents the essential formulations for designing materials with enhanced stiffness and buckling strength using topology optimization.   The finite element method is combined with   homogenization theory used to evaluate material properties~\cite{Cook2002}.  It is well-known that linear elements, i.e., 8-node hexahedral element ($H_8$), overestimate material stiffness and suffer from shear locking, which results in an inaccurate representation of stresses.  To represent the stress situation   more accurately, we employ the incompatible elements proposed by Wilson et al. \cite{Wilson1973}, i.e.,   the so-called $H_{11}$ element.   Three additional so-called incompatible modes are considered to represent bending deformations accurately.  Detailed performance comparison between linear and incompatible elements has been performed and discussed for 2D structures in~\cite{Ferrari2019}.  The reader is referred to Refs~\citep{Wilson1973,Wilson1990} for additional formulations. 
  
 \subsection{Stiffness and buckling strength evaluation}
\begin{figure}[!ht]
 	\centering
  \includegraphics[]{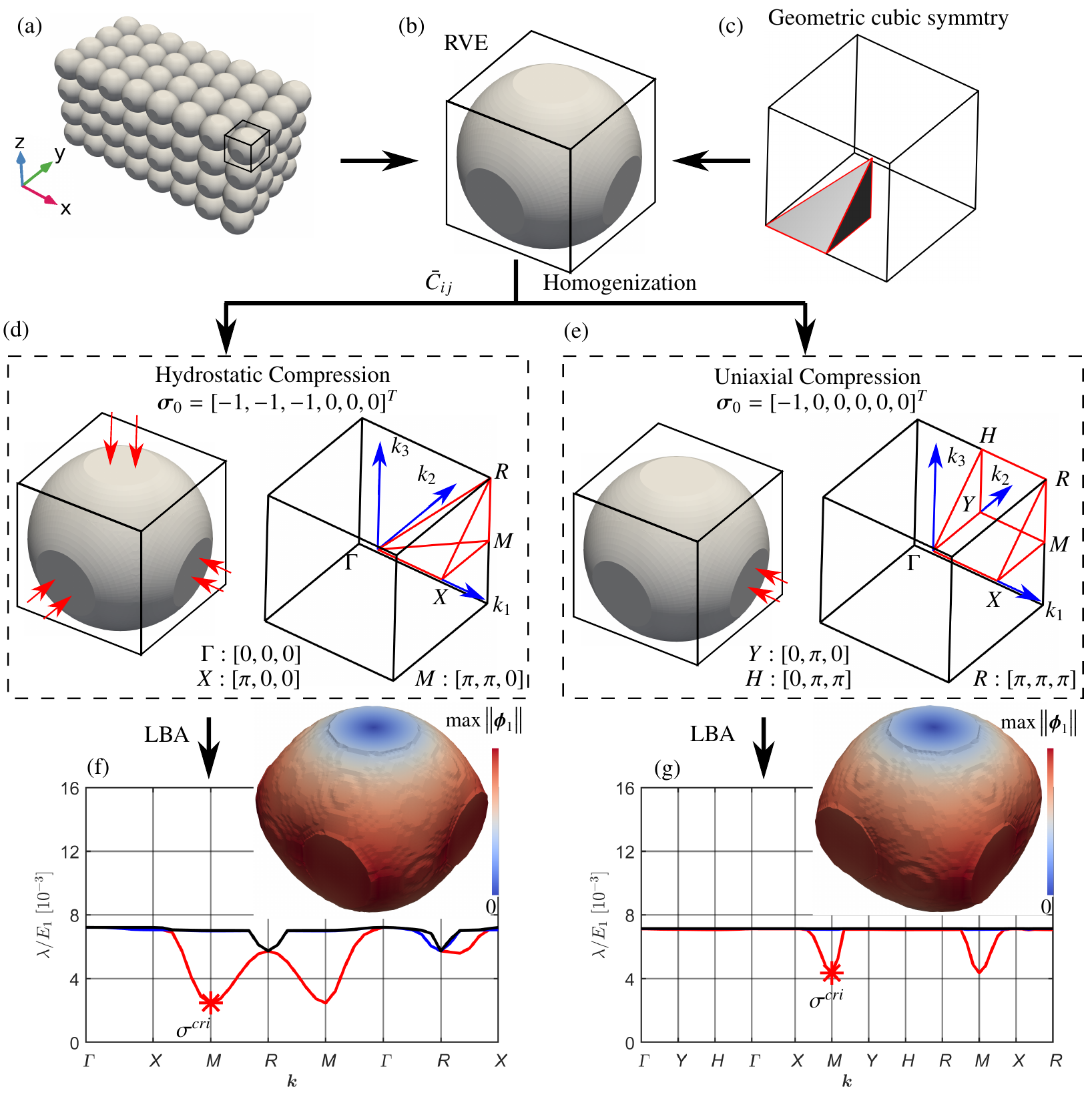}
 	\caption{Schematic Illustration of the modeling procedure and case studies (a)  Periodic material. (b) RVE, i.e., periodic microstructure. (c) Illustration of  geometric cubic symmetry. The repetitive domain in the microstructure: the region enclosed by the red lines. (d)  Illustration of the macroscopic hydrostatic compression and  the corresponding IBZ, i.e., the region enclosed by the red lines.   (e)  Illustration of the macroscopic uniaxial compression and  the corresponding IBZ, i.e., the region enclosed by the red lines.   (f)-(g) Band diagrams of the hollow sphere  under  hydrostatic (f) and uniaxial (g) stress state. Red, blue and black lines represent the first, second and third bands, respectively.   The inserts in (f) and (g) illustrate the corresponding critical buckling modes (inserts). The critical buckling is indicated by   an asterisk  ($\ast$) in the band diagram. The color in the buckling modes indicates the deformation magnitude.  } 	\label{fig:Illustration}
 \end{figure}

 In the small strain limit, the effective material properties can be evaluated using the homogenization approach ~\cite{Sigmund1994,Hassani1998a} based on  a representative volume element (RVE), i.e., the periodic microstructure, see Fig.~\ref{fig:Illustration}. The microstructure dimension is set to 1.  The symmetry properties of the effective elasticity matrix  are exploited to make use of the abbreviation $kl\rightarrow \alpha$: $11\rightarrow 1$, $22\rightarrow 2$, $33\rightarrow 3$,  $(23,32)\rightarrow 4$,  $(13,31)\rightarrow 5$, $(12,21)\rightarrow 6$ to represent the equations in more compact form.   
 
 The  effective elasticity matrix is calculated by an equivalent energy based homogenization formulation~\cite{Sigmund1994,Hassani1998a}:
 \begin{align}
 \bar{C}_{ \alpha\beta} = \frac{1}{|Y|} \sum_{e=1}^N   \int_{Y^e} \left(\tilde{\boldsymbol{\varepsilon}}_{\alpha} - \mathbf{B}^e\boldsymbol{\chi}^e_{\alpha} \right)^T \mathbf{C}^e \left(\tilde{\boldsymbol{\varepsilon}}_{\beta} - \mathbf{B}^e\boldsymbol{\chi}^e_{\beta} \right) dY,
 \end{align}
 where $|Y|$ denotes the volume of the microstructure, $\sum$ represents a finite element assembly operation over all $N$ elements,  the superscript $ \left( \right) ^T$ denotes the transpose,  $\mathbf{B}^e$ with a size of $6\times 24$ is the condensed strain-displacement matrix of element $e$ in the $H_{11}$ element formualtion,  $\mathbf{C}^e$ is the elasticity matrix of the material in element  $e$,  $\tilde{\boldsymbol{\varepsilon}}_{\alpha}=\left[ \delta_{\alpha\beta}\right] $ denotes the 6 independent unit strain fields, and $\boldsymbol{\chi}^e_{\alpha}$ is the perturbation field induced by the $\alpha$th unit strain field under periodic boundaries excluding the incompatible modes, solved by
 \begin{align}\label{eq:chi_disc}
 \mathbf{K}_0 \boldsymbol{\chi}_{\alpha} = \mathbf{f}_{\alpha}, \quad \alpha=1,2,3,4,5,6,\\
 \boldsymbol{\chi}_{\alpha}|_{x=1}=\boldsymbol{\chi}_{\alpha}|_{x=0},   \quad   \boldsymbol{\chi}_{\alpha}|_{y=1}=\boldsymbol{\chi}_{\alpha}|_{y=0},\quad   \boldsymbol{\chi}_{\alpha}|_{z=1}=\boldsymbol{\chi}_{\alpha}|_{z=0}  \nonumber
 \end{align}
Here $\mathbf{K}_0$  is the global  condensed elastic stiffness matrix and   $\mathbf{f}_{\alpha}$ is the condensed equivalent  load vector induced by the $\alpha$th unit strain. The detailed formulation of $ \mathbf{K}_0$  and $\mathbf{f}_{\alpha}$  can be found  in \cite{Thomsen2018,Wilson1990}. 

 The effective material compliance matrix is $\bar{\mathbf{S}}=\bar{\mathbf{C}}^{-1} $. This study only considers microstructures with cubic symmetry that  perform identically  along all the three axes.  However,  everything is   immediately  applicable to anisotropic cases albeit with increased computational efforts.    The effective Young's and bulk moduli and the  isotropy index, i.e., Zener ratio~\cite{Ranganathan2008} are  calculated using  the effective elasticity or compliance matrices, stated as 
 \begin{align}\label{eq:prop}
 \bar{E}=\frac{1}{\bar{S}_{11}},  \quad   
 \bar{\kappa} =\frac{  \bar{C}_{11}+2\bar{C}_{12}}{3}, \quad  a_r=   \frac{2 \bar{C}_{66}
 }{ \bar{C}_{11}- \bar{C}_{12}}.  
 \end{align} 
 
  Subsequently, under a given macroscopic stress state, material buckling strength can be evaluated using LBA with Bloch-Floquet boundary conditions to cover all the possible buckling modes in the material microstructure. For a prescribed  macroscopic stress state, $ \boldsymbol{\sigma}_0$, the effective elasticity matrix is used to transform the macroscopic stress to the macroscopic strain, written as
 \begin{align}
 \boldsymbol{\varepsilon}_0=   \bar{\boldsymbol{S}} \boldsymbol{\sigma}_0.
 \label{eq:EqStrain}
 \end{align}
 Here $\boldsymbol{\sigma}_0=\sigma_0 \boldsymbol{n}$ with $\boldsymbol{n}$ indicating the loading cases and $\sigma_0$ being the stress value, $\boldsymbol{n}=[-1,-1,-1,0,0,0]^T$ for the hydrostatic loading case and $\boldsymbol{n}=[-1,0,0,0,0,0]^T$ for the uniaxial  loading case.   The corresponding stress state of element $e$  is obtained by the superposition of the six perturbation fields induced by the  six unit strain fields, see Eq.~\eqref{eq:chi_disc}, expressed as   
 \begin{align}
 \boldsymbol{\sigma}^e =
 \mathbf{C}^e \Big( \mathbf{I} - \mathbf{B}^e 
 \boldsymbol{X}^e
 \Big)   \boldsymbol{\varepsilon}_0,
 \label{eq:sig0_disc}
 \end{align} 
 where $\boldsymbol{X}^e= \left[\boldsymbol{\chi}^e_{1}\ \ \boldsymbol{\chi}^e_{2}\ \   \boldsymbol{\chi}^e_{3}\ \  \boldsymbol{\chi}^e_{4}\ \   \boldsymbol{\chi}^e_{5}\ \  
 \boldsymbol{\chi}^e_{6} \right] $ is a $24\times 6$ matrix containing the element perturbation fields.
 
 Based on the stress distribution in the microstructure,   subsequent LBA is performed  to evaluate the material buckling strength. Both short  and long wavelength buckling is captured by employing the Floquet-Bloch boundary conditions in the LBA~\cite{Triantafyllidis1993,Neves2002,Thomsen2018}, calculated by
 \begin{align}\label{eq:buckling}
 \left[  \boldsymbol{K}_0 +   \lambda_h  \boldsymbol{K}_{\sigma}  \right]\boldsymbol{\phi}_h = \boldsymbol{0}, \\
 \boldsymbol{\phi}_h|_{x=1}=e^{Ik_1}\boldsymbol{\phi}_h|_{x=0},\quad  \boldsymbol{\phi}_h|_{y=1}=e^{Ik_2}\boldsymbol{\phi}_h|_{y=0},\quad  \boldsymbol{\phi}_h|_{z=1}=e^{Ik_3}\boldsymbol{\phi}_h|_{z=0}. \nonumber
 \end{align}	
Here  $\boldsymbol{K}_{\sigma}=\sum_e \boldsymbol{K}^e_{\sigma} $  is the stress stiffness matrix with $  \boldsymbol{K}^e_{\sigma} $ being the elemental stress stiffness matrix, $I=\sqrt{-1}$ is the imaginary unit, the smallest eigenvalue, $\lambda_1$, is the critical buckling strength for the given wave vector  ($\boldsymbol{k}=[k_1,k_2,k_3]^T$),  and $\boldsymbol{\phi}_1$ is the associated eigenvector. 

From a computational   point of view, it is more convenient to reformulate buckling analysis to    
\begin{equation}\label{eq:buckling1}
\left[ -\tau_h \boldsymbol{K}_0 -   \boldsymbol{K}_{\sigma}  \right]\boldsymbol{\phi}_h = \boldsymbol{0}  
\end{equation}
where $\tau_h=1/\lambda_h$ can be obtained by solving for maximum eigenvalues in the above equation. 
 
 The material buckling strength, $\sigma^{cri}$, is determined by the smallest eigenvalue for all the possible wave vectors in Eq.~\eqref{eq:buckling}, located in the first  Brillouin zone, $\lambda_{min}$, which corresponds to the largest eigenvalue in Eq.~\eqref{eq:buckling1}, i.e., $\sigma^{cri}=\lambda_{min}\sigma_0=1/\tau_{max} \sigma_0 $.  The critical buckling mode is defined as the eigenvector associated with $\tau_{max}$. The first  Brillouin zone is  the primitive cell in reciprocal space~\cite{Brillouin1953}, spanning over  $k_j \in [-\pi,\ \pi], \ j=1,2,3$. It can be further reduced to the IBZ depending on the symmetries shared between the microstructure geometry and the macroscopic stress state.  Fig.~\ref{fig:Illustration} (d) and (e) illustrate the IBZ  for microstructures with cubic symmetry under hydrostatic ($\sigma_0=1$, $\boldsymbol{n}=\left[-1,-1,-1,0,0,0 \right]^T$ ) and uniaxial ($\sigma_0=1$, $\boldsymbol{n}=\left[-1,0,0,0,0,0 \right]^T$) compression, respectively.  The IBZ is denoted  by  the region enclosed by the red lines.  Previous studies have shown that the critical buckling mode can be captured by sweeping $k$-vectors along the boundaries of the IBZ~\cite{Geymonat1993,Thomsen2018,Andersen2020}.  The band diagrams, as customary, are calculated along the boundaries of the IBZ. The critical buckling value is highlighted  by an asterisk  ($\ast  $) in the band diagrams.
 
 The modeling procedure is sketched in~Fig.~\ref{fig:Illustration}.  As a test case, the band diagrams of a hollow sphere are calculated for hydrostatic and uniaxial compressions. The hollow sphere's inner and outer radii are 0.474 and 0.541, and the corresponding volume fraction is 0.2.  The properties of the base material are $E_1=1$ and $\nu=1/3$.  The effective microstructure properties including Young's,  bulk moduli and buckling strength are presented as relative quantities with respect to  the base material Young's  modulus, i.e., $\bar{E}/E_1$, $\bar{\kappa}/E_1$ and $\sigma^{cri}/E_1$, respectively.    The microstructure is discretized by $64 \times 64 \times 64$ cubic $H_{11}$ elements.  The effective properties of the hollow sphere are found as $ \bar{E}/E_1= 0.0316 $, $ \bar{\kappa}/E_1=0.0201  $ and  $a_r=1.23$.   Fig.~\ref{fig:Illustration} (f) and (g) present the corresponding band diagrams under hydrostatic and uniaxial compression, respectively.   It is seen that critical buckling in both cases occurs at the M point of  $\boldsymbol{k}=\left[\pi,\pi,0 \right]$, corresponding to anti-periodic in the x- and y-direction and periodic in the z-direction.   The critical buckling modes (see the inserts) also illustrate that the hollow sphere buckles in similar mode shapes in both loading conditions. The color in the buckling modes indicates the deformation magnitude unless otherwise stated. The buckling strength under uniaxial compression ($\sigma^{cri}/E_1=0.00436$) is higher than the buckling strength under hydrostatic compression ($\sigma^{cri}/E_1= 0.00247$).

  \subsection{Design parameterization}
An element-wise constant set of physical design variables, $\bar{\rho}_e$, is employed to represent the material distribution in the microstructure.  Element  $e$  is occupied by the base material when $\bar{\rho}_e=1$ and by void when  $\bar{\rho}_e=0$.  As in \cite{Thomsen2018}, to suppress the spurious buckling modes associated with the low stiffness elements, different interpolation schemes are employed for the elastic stiffness and stress stiffness using the solid isotropic material with penalization (SIMP) model \citep{Bendsoe1999}, stated as:
\begin{equation} \label{Eq:2DCon}
\mathbf{E}_e = 
\left\{
\begin{array}{ll}
\bar{\rho}^p_e(E_1-E_0)+E_0 & \text{for } \boldsymbol{K}_0,   \mathbf{f}_{\alpha} \\
\bar{\rho}^p_e  \mathbf{E}_1 & \text{for } \boldsymbol{K}_{\sigma}
\end{array}
\right.
\end{equation} 
where  $E_0=10^{-4} E_1$ is employed to represent void regions,   $p=3$ is chosen as the standard penalization factor.  

In the density-based topology optimization approach,  an element-wise constant set of design variables ${\rho}_e$ is introduced to generate the  material distribution in the microstructure.   A hyperbolic tangent threshold projection, is employed to generate physical design variables from the design variables,  ${\rho}_e$, to enhance  the discreteness  of the optimized design \citep{Wang2011}. This is given
\begin{equation} \label{Eq:2DProj}
\bar{\rho}_e=\frac{\tanh{\left(\beta_1\eta\right)}+\tanh{\left(\beta_1\left(\widetilde{\rho}_e-\eta\right)\right)}}{\tanh{\left(\beta_1\eta\right)}+\tanh{\left(\beta_1\left(1-\eta\right)\right)}}
\end{equation}
where $\widetilde{\rho}_e$ is the filtered design variable calculated from the design variables, $\boldsymbol{\rho}$, using a PDE filter presented in \citep{Lazarov2011b}. The reader is referred to~\cite{Lazarov2011b} for further details.   When $\beta_1$ is big,   $\bar{\rho}_e\approx 1$ if $\widetilde{\rho}_e>\eta$  and $\bar{\rho}_e\approx0$ if $\widetilde{\rho}_e<\eta$. Hence the projection in Eq. \eqref{Eq:2DProj}  suppresses  gray element density regions induced by the PDE filter when  $\beta_1$ is big and ensures black-white designs when the optimization converges.  Moreover, it mimics the manufacturing process, and is used in the robust design formulations context~\citep{Wang2011}, where the manufacturing errors are taken into accounts by choosing different thresholds, $\eta$, as discussed later.

\subsection{Design problem formulation}

The Kreisselmeier–Steinhauser (KS) function~\cite{Kreisselmeier1980} is employed  to aggregate the considered eigenvalues for  given $\boldsymbol{k}$-vectors to represent the material buckling strength as
\begin{align}\label{eq:KS}
KS (\tau_h\left(\boldsymbol{k}_l \right) ) = \frac{1}{\zeta} \ln\left(  \sum_{l=1}^{n_h} \sum_{h=1}^{m_l}e^{\zeta \tau_h\left( \boldsymbol{k}_l\right)}\right) .   
\end{align}
The optimization problem for enhancing material stiffness and  buckling strength can be formulated to minimize a weighted value of the KS function and $\bar{E}^{-1}$, stated as
\begin{align}
\min\limits_{\boldsymbol{\rho}} \quad   &  \qquad  \gamma_1	KS  \left( \tau_h\left(\boldsymbol{k}_l \right)  \right) + \left( 1-\gamma_1 \right) \bar{E}^{-1}   \nonumber \\
s.t. & \qquad  \left[ -\tau_h  \boldsymbol{K}_0 \left(\boldsymbol{k}_l \right) -    \boldsymbol{K}_{\sigma} \left( \boldsymbol{k}_l\right)  \right]\boldsymbol{\phi}_h = \boldsymbol{0} \nonumber \\
& \qquad \mathbf{K}_0 \boldsymbol{\chi}_{ \alpha} = \mathbf{f}_{\alpha}, \quad \alpha=1,2,3,4,5,6   \\ 
& \qquad   \left(  a_r-1\right)^2<\delta^2    \nonumber\\
&  \qquad f=\frac{\boldsymbol{v}^T\bar{\boldsymbol{\rho}} \left(\boldsymbol{\rho}\right)}{\sum\limits_e v_e}   \leq f^*  \nonumber \\
&  \qquad  \boldsymbol{0} \leq \boldsymbol{\rho} \leq \boldsymbol{1} \nonumber
\label{eq:opt}
\end{align}
where  $\gamma_1$ is the weight, $\delta$ is the allowed error on the material isotropy,   $\boldsymbol{v}$ is the elemental volume vector with $v_e$ representing the volume of element $e$, and  $f$ and $f^*$ are the actual volume fraction and    the prescribed upper bound of the volume fraction in the microstructure. 

The sensitivity of a component, $\bar{C}_{ \alpha\beta}$, in the  effective elastic matrix $\bar{\boldsymbol{C}}$    with respect to   $\bar{\rho}^e$, is written as
\begin{equation}\label{eq:CCSens}
\frac{\partial \bar{C}_{\alpha\beta} } {\partial \bar{\rho}^e} = \frac{1}{|Y|}     \int_{Y^e}  \left(\tilde{\boldsymbol{\varepsilon}}_{\alpha} - \mathbf{B}^e\boldsymbol{\chi}^e_{\alpha} \right)^T \frac{\partial  \mathbf{C}^e } {\partial \bar{\rho}^e}\left(\tilde{\boldsymbol{\varepsilon}}_{\beta} - \mathbf{B}^e\boldsymbol{\chi}^e_{\beta} \right).  
\end{equation}
The sensitivities of  $\bar{E}$ and $a_r$ can be analytical derived using Eqs. \eqref{eq:prop} and  \eqref{eq:CCSens}.

Assuming  a  distinct eigenvalue, $\tau_h$, and that the eigenvector is normalized, as $\left( \boldsymbol{\phi}_h \right)^H \boldsymbol{K}_0   \boldsymbol{\phi}_h = 1  $, the sensitivity of eigenvalue $\tau_h$ with respect to  $\bar{\rho}^e$ can be obtained via the adjoint sensitivity analysis as described below~\cite{Thomsen2018}, 
\begin{align}
\frac{\partial \tau_h}{\partial \bar{\rho}^e}=&
{\left( \boldsymbol{\phi}^e_h\right)}^H  \left[ -\tau_h \frac{\partial \boldsymbol{K}^e_{0}  }{\partial \bar{\rho}^e} -  \frac{\partial \boldsymbol{K}^e_{\sigma} } {\partial \bar{\rho}^e} \right]  \boldsymbol{\phi}^e_h + {\left( \boldsymbol{\phi}_h\right)}^H \left[  -  \frac{\partial \boldsymbol{K}_{\sigma} }{ \partial  \boldsymbol{\varepsilon}_0} \frac{ \partial  \boldsymbol{\varepsilon}_0}{\partial \bar{\rho}^e} \right] \boldsymbol{\phi}_h \\
&+\sum_{\alpha=1}^6\left( \boldsymbol{\psi}_ \alpha^e\right)^H \left[ 
\frac{\partial \boldsymbol{K}^e_{0}  }{\partial \bar{\rho}^e} \boldsymbol{\chi}^e_{ \alpha} - \frac{\partial \mathbf{f}^e_{\alpha}}{\partial \bar{\rho}^e} \right]. 
\label{eq:Sens}
\end{align} 
Here $\left( \right) ^H$ denotes the complex conjugate,  ${ \partial  \boldsymbol{\varepsilon}_0}/{\partial \bar{\rho}^e}$ can be directly derived using \eqref{eq:EqStrain} and \eqref{eq:CCSens},  and  $\boldsymbol{\psi}_ \alpha$ is the adjoint vector corresponding to $\boldsymbol{\chi}_{\alpha}$, which is obtained by
\begin{align}
\boldsymbol{K}_{0}    \boldsymbol{\psi}_\alpha=\sum_e  \left(\boldsymbol{\phi}^e_h\right)^H   \left[  \frac{\partial \boldsymbol{K}^e_{\sigma} }{ \partial  \boldsymbol{\chi}^e_\alpha}  \right] \boldsymbol{\phi}^e_h 
\end{align} 
The reader is referred to \cite{Thomsen2018} for the detailed calculation of $ \frac{\partial \boldsymbol{K}^e_{\sigma} }{ \partial  \boldsymbol{\chi}^e_ \alpha}$.

One of the advantages of using aggregation functions is the uniqueness of the eigenvalue gradient, even for repeated eigenvalues, as long as eigenvectors are   $\boldsymbol{K}_0$-orthonormalized. The gradient uniqueness has been proved for the p-norm aggregation function in~\cite{Torii2017}. The same argument can be trivially extended to the KS-function used in this study.  The sensitivities of the objective and constraints with  respect to  design variable, ${\rho}_e$,  are obtained using the chain rule.

To enhance design robustness and realize length scale control on the  microstructure, the one case robust formulation is employed in this study as in \cite{Lazarov2016}. The objective is evaluated on an eroded microstructure generated using $0.5+\Delta \eta$ with the volume constraint working on the dilated  microstructure  generated with a threshold of $0.5- \Delta\eta$. The volume constraint is updated every 20 iterations such that the volume constraint of the intermediate microstructure of $\eta=0.5$ is satisfied.  In order to ensure black-white designs, a continuation scheme is employed  to increase $\beta_1$ in Eq. \eqref{Eq:2DProj}.  $\beta_1$ is updated every 40 iterations using $\beta_1=2  \beta_1$, with a maximal value of $\beta_{\rm{max}}=50$ and initial value of $\beta_1=1$. The number of eigenvalues considered for a given $\boldsymbol{k}_l$, $m_l$, is set to the number of eigenvalues within 5\% of the maximum $\tau_h$ for all $\boldsymbol{k}$-vectors considered. The number of eigenvalues   is updated every 20 iterations. The exponent in Eq.~\eqref{eq:KS} is set to $\zeta=100/\max\left(\tau_h\left( \boldsymbol{k}_l\right) \right) $ and updated every 100 iterations. 

The optimization problem is implemented in  a flexible framework for large scale topology optimization  \cite{Aage2015} using the Portable Extensible Toolkit for Scientific computation (PETSc)~\cite{Balay2016} and Scalable Library for Eigenvalue Problem Computations  (SLEPc)~\cite{Hernandez2005}. The design is iteratively updated using the Method of Moving Asymptotes (MMA) \cite{Svanberg1987} based on the gradients of the objective and constraints.

\section{Results}
The proposed optimization formulation is employed to design 3D isotropic microstructures with enhanced stiffness and buckling strength.  The unit cell is discretized by $64\times 64 \times  64$ $H_{11}$ elements.  The one-case robust formulation is employed with  $\Delta \eta=0.05$  and a filter radius of $r=0.05$, which corresponds to a relative minimal length scale around 0.02, i.e., the minimal feature size is around 2\% of the microstructure size~\cite{Silva2021}.  The volume fraction upper bound on the intermediate microstructure is set to $f^*=0.2$, and the allowed error on the Zener ratio is $\delta=0.05$.   The initial design is the hollow sphere shown in Fig.~\ref{fig:Illustration}.  All the optimized microstructures including shape optimized ones are obtained with the grayness less than $1\%$, which is defined as $ \sum_e  4 \left(1-\bar{\rho}^e \right) \bar{\rho}^e /\sum_e   1 $ as in~\cite{Wang2011}, and they are shown as smoothed configurations from a post-process using an iso-volume threshold of 0.5 on the filtered density  in  Paraview~\cite{Ahrens2005}, unless otherwise stated.

\subsection{Microstructure with enhanced buckling strength under hydrostatic compression}
In the first case, we design isotropic microstructure with enhanced buckling strength ($\gamma_1=1$) under hydrostatic compression of   $\boldsymbol{\sigma}_0=[-1,-1,-1,0,0,0]^T$.  Five $\boldsymbol{k}$-vectors are chosen as the target points in the IBZ, i.e., four vertices and    $\boldsymbol{k}=[\pi/20,0,0]$ for capturing global shear buckling modes.

 \begin{table}[!htb]
	\centering
	 \includegraphics[]{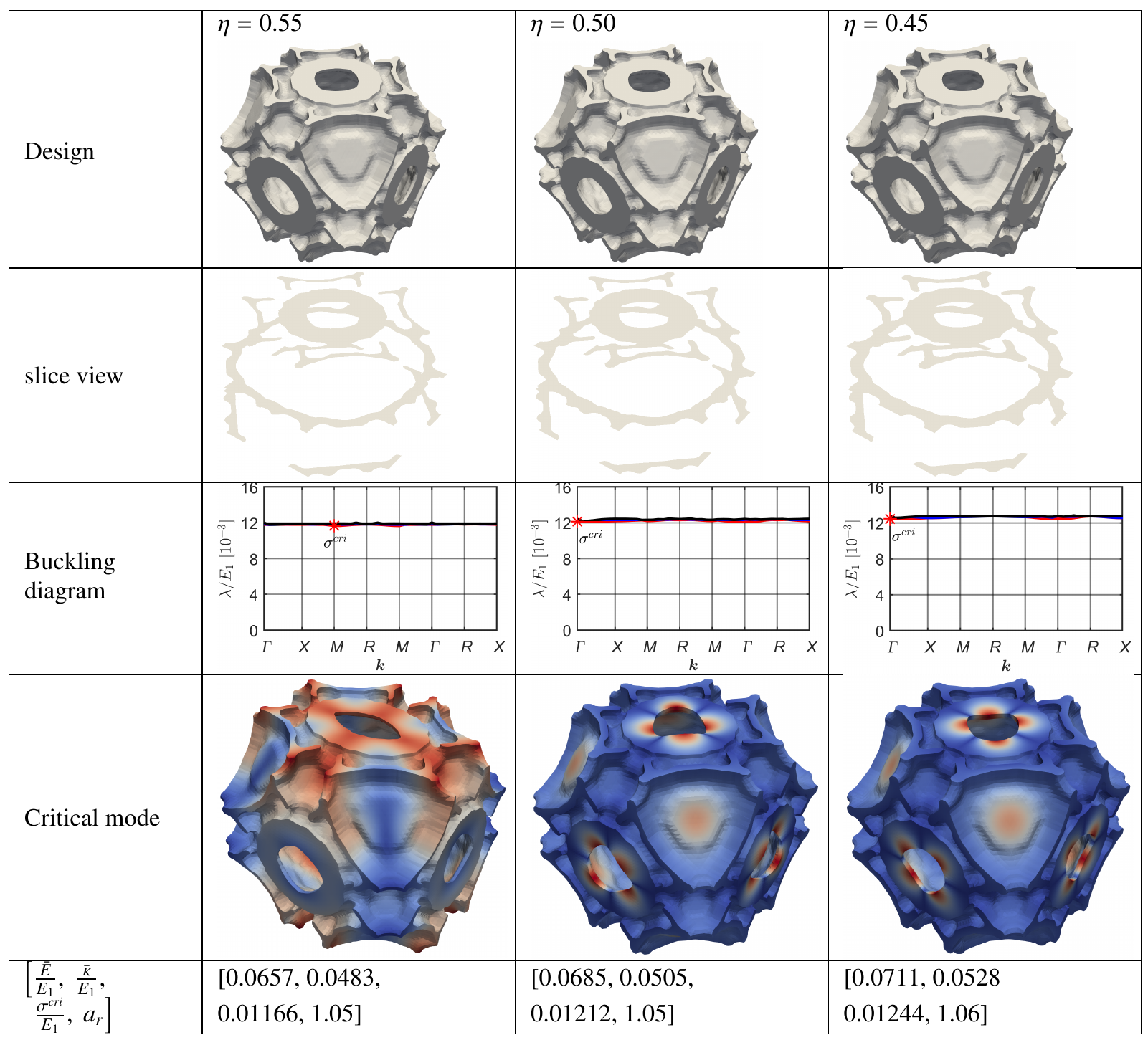}
	\caption{Optimized microstructure with enhanced buckling strength under hydrostatic compression   and corresponding performances. Column 2: Eroded microstructure, Column 3: Intermediate microstructure, Column 4: Dilated microstructure. The sphere microstructure in Fig.~\ref{fig:Illustration} was used as starting guess} \label{tab:Designhydro}
\end{table}

Tab.~\ref{tab:Designhydro} summarizes the performances of the three design realizations considered in the robust optimization process.  It is seen that the three design realizations share the same topology. Hence the intermediate design possesses a length scale in both solid and void regions \cite{Wang2011}.  The solid region's thickness increases from the eroded to the dilated designs, which results in increasing moduli and buckling strength (see the last row in Tab.~\ref{tab:Designhydro}).   The slice views in the second row in Tab.~\ref{tab:Designhydro} illustrate that the optimized  microstructure possesses a single length-scale. In contrast to the study~\cite{Thomsen2018},  where the optimizer exploits geometric hierarchy to improve material buckling strength,  the development of geometric hierarchy is here prevented by the imposed length scale for the given volume fraction $20\%$.  

The first three buckling bands along the IBZ  boundaries exhibit minor variations and are tightly clustered. It shows that buckling for short- and long-wavelengths in different directions coincides, implying that the optimized microstructure has made near-optimal use of available material to equally strengthen the different buckling modes.  In contrast, the initial design's buckling bands display significant variation, indicating that a single failure mode is dominating and that the material is non-optimal in terms of buckling strength.  It is noted that the intermediate and dilated microstructures fail due to a cell-periodic buckling mode, i.e., at the $\Gamma$ point,  while the eroded microstructure fails at the M point, i.e., an x-, y-antiperiodic and z-periodic mode.   Due to the minor variations in the first buckling band,   the proposed one-case robust optimization works as expected, i.e., the intermediate microstructure's buckling strength has been significantly improved by the proposed optimization.  The optimized microstructure exhibits an improvement of $543\%$  in buckling strength and $125\%$ and  $162\%$  in Young's and bulk moduli compared to the starting guess, respectively. The detailed comparisons demonstrate the enhancement of the buckling strength is achieved by splitting a large region into several small regions and opening the closed-wall cell.  These geometrical changes enlarge the thickness-to-area ratio, which increases effective plate bending  stiffness and hence results in superior buckling strength.

\subsection{Microstructures with enhanced stiffness and  buckling strength under uniaxial compression}

In the second case,  we  design isotropic microstructures with enhanced stiffness and  buckling strength under uniaxial compression    $\boldsymbol{\sigma}_0=[-1,0,0,0,0,0]^T$.  Seven $\boldsymbol{k}$-vectors are chosen as the target points here, i.e., six vertices in the IBZ and  $\boldsymbol{k}=[\pi/20,0,0]$ for capturing global shear buckling modes. 

\begin{figure}[!htb]
	\centering
	\includegraphics[]{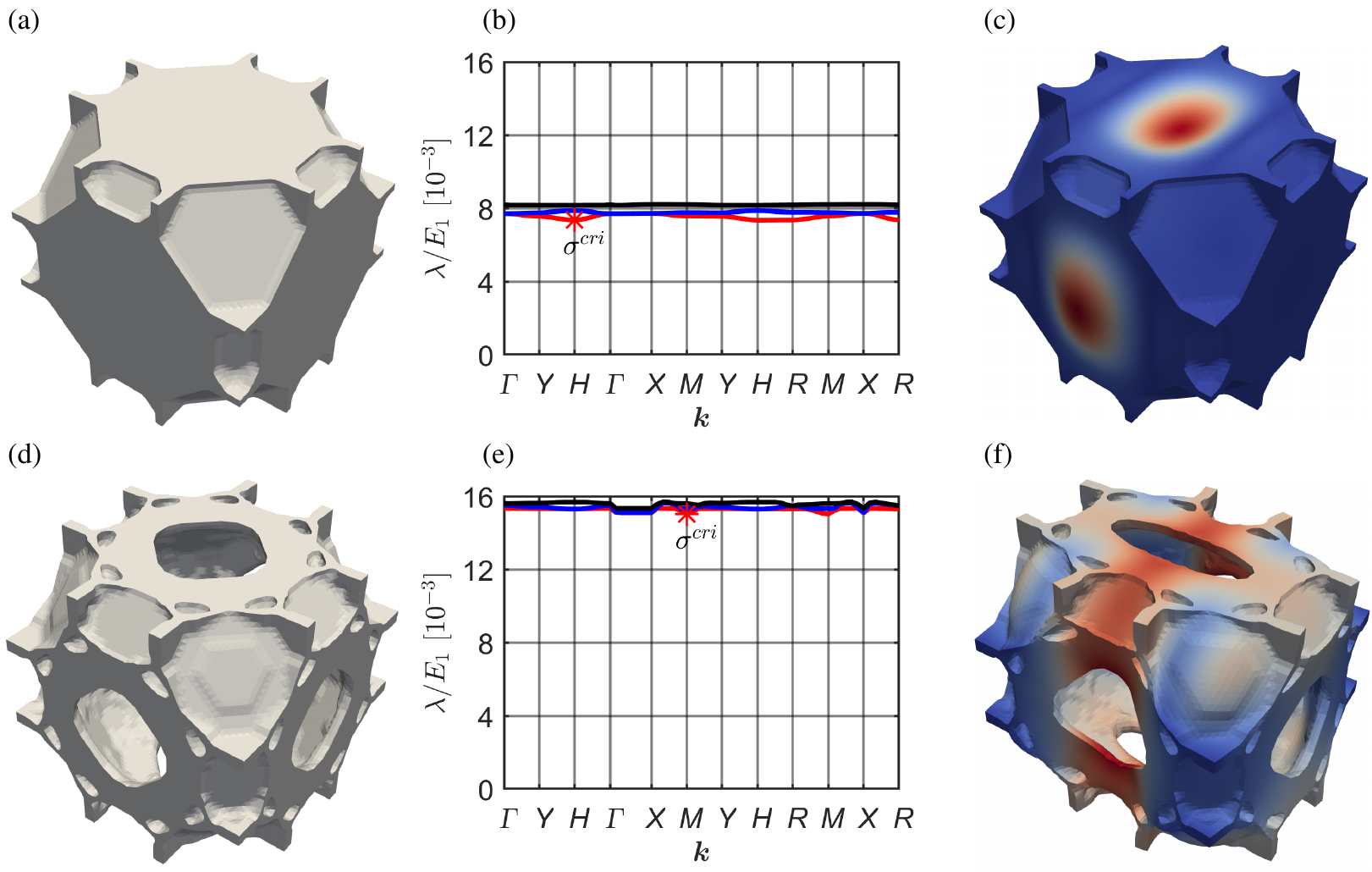}
	\caption{The performance of the optimized  microstructures for maximum  Young's modulus $\bar{E}$    (a)-(c) and for maximum buckling strength  $\sigma^{cri}$  under uniaxial compression (d)-(f). Left column: Optimized microstructures and effective Young's moduli of the two microstructures are $\bar{E}/E_1=0.0884$ (a) and $\bar{E}/E_1=0.0678$ (d); Middle column: Band diagram; Right column: Critical buckling mode. }
	\label{fig:OptDes}
\end{figure}

Fig.~\ref{fig:OptDes}  presents the optimized microstructures  for  maximum Young's modulus  $\bar{E}$ ($\gamma_1=0$  ), see (a)-(c), and the one for  maximum buckling strength $\sigma^{cri}$  ($\gamma_1=1$  ) see (d)-(f). The microstructure with the maximum   $\bar{E}$  is a closed-wall cell and consists of flat plates.  It is seen that the critical buckling mode is located at the H point with  $\boldsymbol{k}=[0,\pi,\pi]$, corresponding to   anti-periodic in the y- and z- direction and periodic in the x-direction. The corresponding buckling mode in (c) demonstrates that the material's buckling strength is determined by the plates' buckling normal to the $y$- or $z$-axis.   In contrast to the material with the maximum  $\bar{E}$, the optimized microstructure with the maximum  $\sigma^{cri}$ consists of   more complex plates which locally adapt thicknesses that increase the buckling strength.  A small variation is observed in the first buckling band, and the corresponding deviation is within $3\%$. The critical buckling mode is located at the M point with  $\boldsymbol{k}=[\pi,\pi,0]$. However, the buckling modes with a long wavelength occur almost simultaneously along $\Gamma-X$ in the band diagram, see (e).

\begin{figure}[!htb]
	\centering
	\includegraphics[ ]{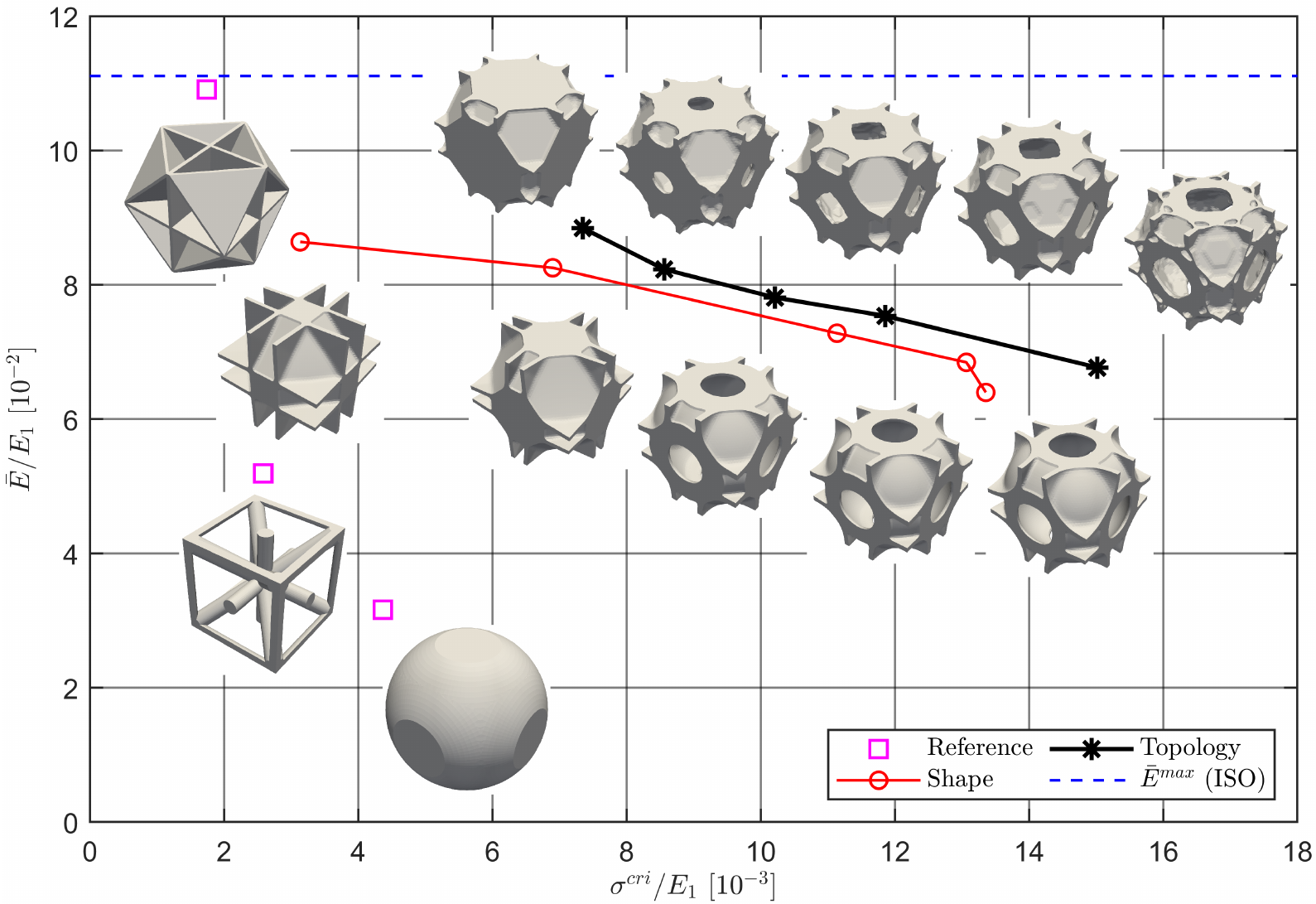}
	\caption{Young's modulus vs. uniaxial buckling strength  response  of the  topology optimized (black line) and shape optimized (red line) microstructures with tunable stiffness and buckling strength and the reference microstructures, i.e., isotropic stiffness-optimal plate and truss microstructures from~\cite{Berger2017} and the initial spherical microstructure. }
	\label{fig:OptDesUniShapeAllE}
\end{figure}

The properties of the optimized microstructure for maximum $\bar{E}$  are $\bar{E}/E_1=0.0884$,   $\bar{\kappa}/E_1=0.0744 $,  ${\sigma^{cri}}/{E_1}=0.00735$.
Compared to the theoretical Hashin–Shtrikman upper bounds of isotropic material with $f=20\%$,  ($\bar{E}^{max}/E_1= 0.1111$  and  $\bar{\kappa}^{max}/E_1=   0.0769$), the optimized material microstructure for  maximum   $\bar{E}$ is   a  bulk-modulus-optimal and sub-Young's-modulus-optimal microstructure although optimized for  maximum   $\bar{E}$. The Young's-modulus-optimal microstructure is hindered by the imposed length scale for  $20\%$ volume fraction, because the stiffness-optimal plate microstructure (see the plate microstructure at the upper left corner in Fig.~\ref{fig:OptDesUniShapeAllE}) has thinner plates than allowed by the mesh and imposed length-scale.

 \begin{figure}[htb!]
 	\centering
 		\includegraphics[ ]{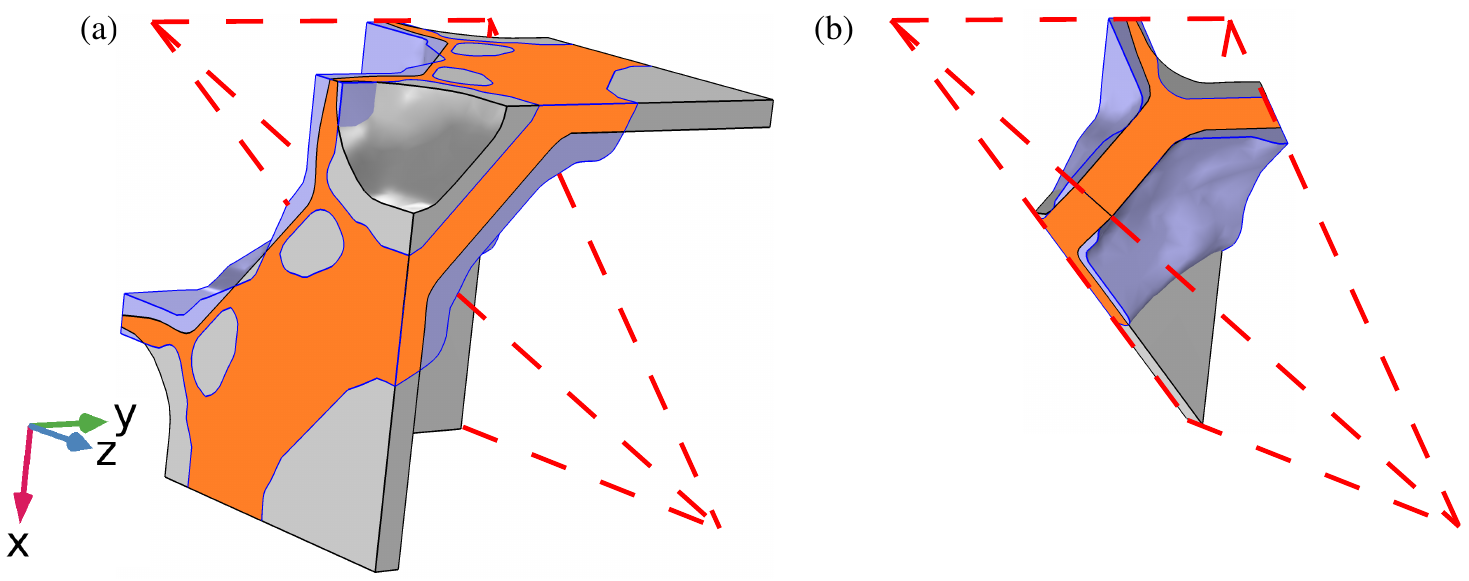}
 	\caption{Geometric comparison  between the  optimized mirostructures for maximum $\bar{E}$ and ${\sigma^{cri}}$ illustrated in 1/8 cell (a) and 1/48 cell (b). Orange: Region  shared by the two optimized microstructures. Purple: Regions owned by the optimized microstructure  for maximum ${\sigma^{cri}}$. Grey: Regions owned by the optimized mirostructure  for   maximum $\bar{E}$. }
 	\label{fig:GeoDiff}
 \end{figure}

The properties of the optimized microstructure for  the maximum  $\sigma^{cri}$ are $\bar{E}/E_1=0.0678$, $\bar{\kappa}/E_1=0.0581 $, $\sigma^{cri}/{E_1}= 0.01503 $. Detailed comparisons between the  two microstructure geometries (see Fig. \ref{fig:GeoDiff}) demonstrate the buckling strength enhancement is attributed to two aspects, 1) thickness increase of some members  by   avoidance of the buckling part in the flat plates;    2)  the evolution of the adapted plates from the flat plates along the body diagonal direction. These geometrical changes lead to the increase of the bending stiffness of the microstructure, hence  result in a relative improvement of $105\%$ in terms of buckling strength compared to the microstructure for maximum $\bar{E}$.  However they also induce  a relative reduction of  $24\%$  in Young's modulus. 

By choosing different weighting factors $\gamma_1$, we design  microstructures with tunable stiffness and buckling strength. The topology optimized microstructures and corresponding performances are summarized in Fig.~\ref{fig:OptDesUniShapeAllE},  together with the shape optimized microstructures discussed in the subsequent subsection, as well as  the initial design, isotropic  stiffness-optimal plate and  truss microstructures from~\cite{Berger2017} with a volume fraction of $20\%$.  The black line with asterisks  summarizes the performance of the topology optimized microstructures. The asterisks represent microstructures obtained by reducing $\gamma_1$ from the right to the left side.  The corresponding  $\gamma_1$ is    0,  0.005,   0.01,  0.02, and 1.   It is remarked that the void areas in the plates normal to the axes get smaller, and  the flatness of the plates also  increases as $\gamma_1$ decreases.  These geometrical changes lead to higher stiffness and lower buckling strength. Compared to the stiffness-optimal plate microstructure from~\cite{Berger2017},  the optimized microstructure for  maximum $\bar{E}$ exhibits $366\%$ improvement in buckling strength with $19\%$ reduction in Young's modulus, and the optimized microstructure for maximum $\sigma^{cri}$  exhibits  $763\%$  improvement in buckling strength with $38\%$ reduction in Young's modulus. The optimized microstructures show better performance in  both Young's modulus and buckling strength than the isotropic  truss microstructure.

\subsection{Simplification of microstructure configurations based  on  shape optimization}
\begin{figure}[!htb]
	\centering
 \includegraphics[]{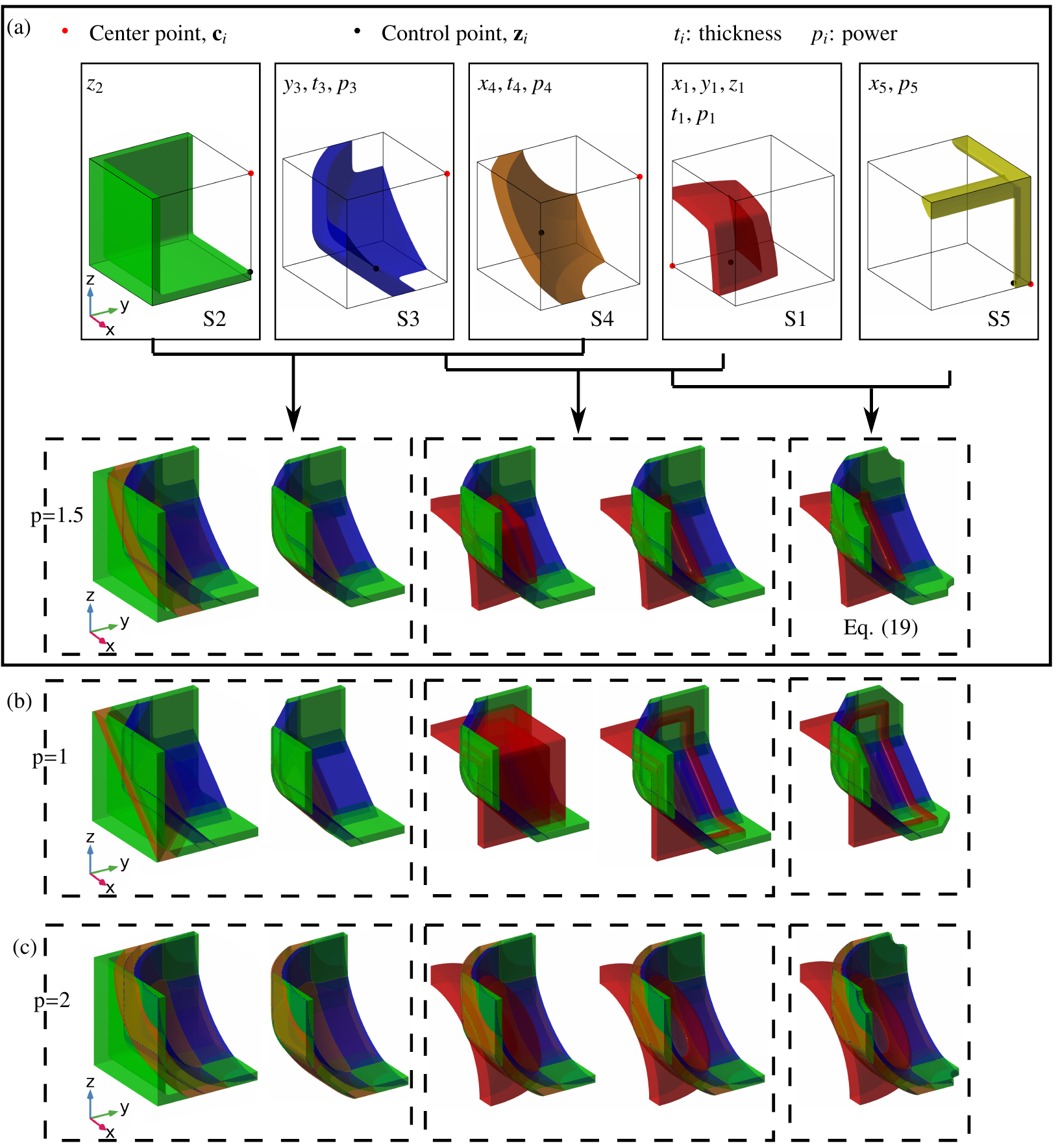} 
	\caption{Illustration of the generation process of  microstructures with different powers using feature-based shape parameterizations.  (a) 
		Illustration of the generation process of the shape parameterized microstructure in Tab.~\ref{Tab:Para}  with $p_i=1.5$ (Second row) together with an illustration of the corresponding four hollow and one solid super-ellipsoids (first row).  (b) Microstructure generation with $p_i=1$. (c)   Microstructure generation with $p_i=2$. Last column: the final microstructures generated by the shape parameterization, i.e., Eq.~\eqref{eq:SupFin} together with cubic symmetry. }
	\label{fig:ShapePara}
\end{figure}

Inspired by the topology optimized microstructure configurations in Fig.~\ref{fig:OptDesUniShapeAllE}, we further simplify the topology optimized microstructure architecture  using a subsequent  feature-based shape optimization.  In the feature-based approach, the physical density variables,  $\bar{\boldsymbol{\rho}}$ are generated via  Boolean operations on a set of hollow  and a solid  super-ellipsoids.  These super-ellipsoids are employed to parameterize 1/48 microstructure and  the rest is mapped according to the cubic symmetry.   As illustrated in the first row  in Fig.~\ref{fig:ShapePara} (a), a hollow super-ellipsoid is controlled by a center point,  $\mathbf{c}_{i}=\left[x_{0i},\ y_{0i},\ z_{0i}\right]$, a control point located at the inner super-ellipsoid  surface which controls the semi-diameters, $\mathbf{z}_{i}=\left[x_{1i},\  y_{1i},\  z_{1i}\right]$, a thickness $t_i$ and the power $p_i$.  We first introduce a function, $T(\mathbf{x},\mathbf{c}_i,\mathbf{z}_{i},p_i,t)$ defined as
 \begin{equation} \label{eq:TPre}
T(\mathbf{x},\mathbf{c}_i,\mathbf{z}_{i},p_i,t)= \sqrt[p_i]{ \left|  {n_{1i}\left( x-x_{0i} \right)}   \right|^{p_i}+ \left|  {n_{2i} \left(y-y_{0i} \right)}  \right| ^{p_i}+ \left|   {n_{3i} \left( z-z_{0i} \right)}\right| ^{p_i}} +t
\end{equation}
where $\mathbf{x}=\left[x,\  y,\  z\right]$ represents any given location, $ \left[  n_{1i},  n_{2i},n_{3i} \right] = \left[ \mathbf{z}_{i}-\mathbf{c}_{i}\right]/d_i$ is the normalized direction vector between the control and center points,  and   $d_i=\parallel \mathbf{z}_{i}- \mathbf{c}_{i}\parallel+\delta_0$ is the distance between the center point and  the control point with $\delta_0=10^{-10}$ to avoid zero denominator.   The hollow super-ellipsoid with thickness $t$ is defined by 


\begin{equation}
H(\mathbf{x},\mathbf{c}_i,\mathbf{z}_{i},p_i,t)= {T(\mathbf{x},\mathbf{c}_i,\mathbf{z}_{i},p_i,0) }-{T(\mathbf{z}_i,\mathbf{c}_i,\mathbf{z}_{i},p_i,t)} \label{HPre}
\end{equation}
Here $T(\mathbf{z}_i,\mathbf{c}_i,\mathbf{z}_{i},p_i,t)$ denotes  the function value at the control point, $\mathbf{z}_{i}$.  $t=0$ and $t=t_i$ in $H(\mathbf{x},\mathbf{c}_i,\mathbf{z}_{i},p_i,t)$  represent the inner and outer super-ellipsoid surfaces, respectively.  The   semi-diameters  of the inner super-ellipsoid along the x-, y- and z-direction are $T(\mathbf{z}_i,\mathbf{c}_i,\mathbf{z}_{i},p_i,0)/n_{1i}$, $T(\mathbf{z}_i,\mathbf{c}_i,\mathbf{z}_{i},p_i,0)/n_{2i}$ and $T(\mathbf{z}_i,\mathbf{c}_i,\mathbf{z}_{i},p_i,0)/n_{3i}$, respectively.  If the location $\mathbf{x}$ is located inside the $i$th inner super-ellipsoid, $T(\mathbf{x},\mathbf{c}_i,\mathbf{z}_{i},p_i,0)<T(\mathbf{z}_i,\mathbf{c}_i,\mathbf{z}_{i},p_i,0)$, hence $H(\mathbf{x},\mathbf{c}_i,\mathbf{z}_{i},p_i,0)<0$, if outside, $T(\mathbf{x},\mathbf{c}_i,\mathbf{z}_{i},p_i,0)>T(\mathbf{z}_i,\mathbf{c}_i,\mathbf{z}_{i},p_i,0)$, hence $H(\mathbf{x},\mathbf{c}_i,\mathbf{z}_{i},p_i,0)>0$ and on  the $i$th super-ellipsoid, $H(\mathbf{x},\mathbf{c}_i,\mathbf{z}_{i},p_i,0)=0$. The densities of the element,  $e$, resulting from the $i$th inner and outer super-ellipsoids  are defined as,
\begin{equation}\label{eq:Sup}
{\rho}^{i,0}_e = \frac{1}{1+\exp\left(\beta_2 H(\mathbf{x}_e,\mathbf{c}_i,\mathbf{z}_{i},p_i,0)\right)}, \qquad
{\rho}^{i,1}_e = \frac{1}{1+\exp\left(\beta_2 H(\mathbf{x}_e,\mathbf{c}_i,\mathbf{z}_{i},p_i,t_i)\right)}
\end{equation}
where $\mathbf{x}_e$ is the centroid of element $e$,  ${\rho}^{i,0}_e$   and   ${\rho}^{i,1}_e$  denote densities generated by the $i$th inner and outer super-ellipsoids, respectively, ${\rho}^{i,0}_e=1$  if $\mathbf{x}_e$ is inside the $i$th inner super-ellipsoid and ${\rho}^{i,0}_e=0$  outside of the $i$th  inner super-ellipsoid.   ${\rho}^{i,1}_e=1$  if $\mathbf{x}_e$ is inside the $i$th outer super-ellipsoid and ${\rho}^{i,1}_e=0$  outside of the $i$th  outer super-ellipsoid, and $\beta_2=200$ is the projection value.

\begin{table}[!htb]
	\centering
	\begin{tabular}{m{1cm} m{2.5cm} m{2.5cm}  m{2.5cm} m{2.5cm} m{2.5cm}  }
	\multicolumn{6}{c}{\includegraphics[]{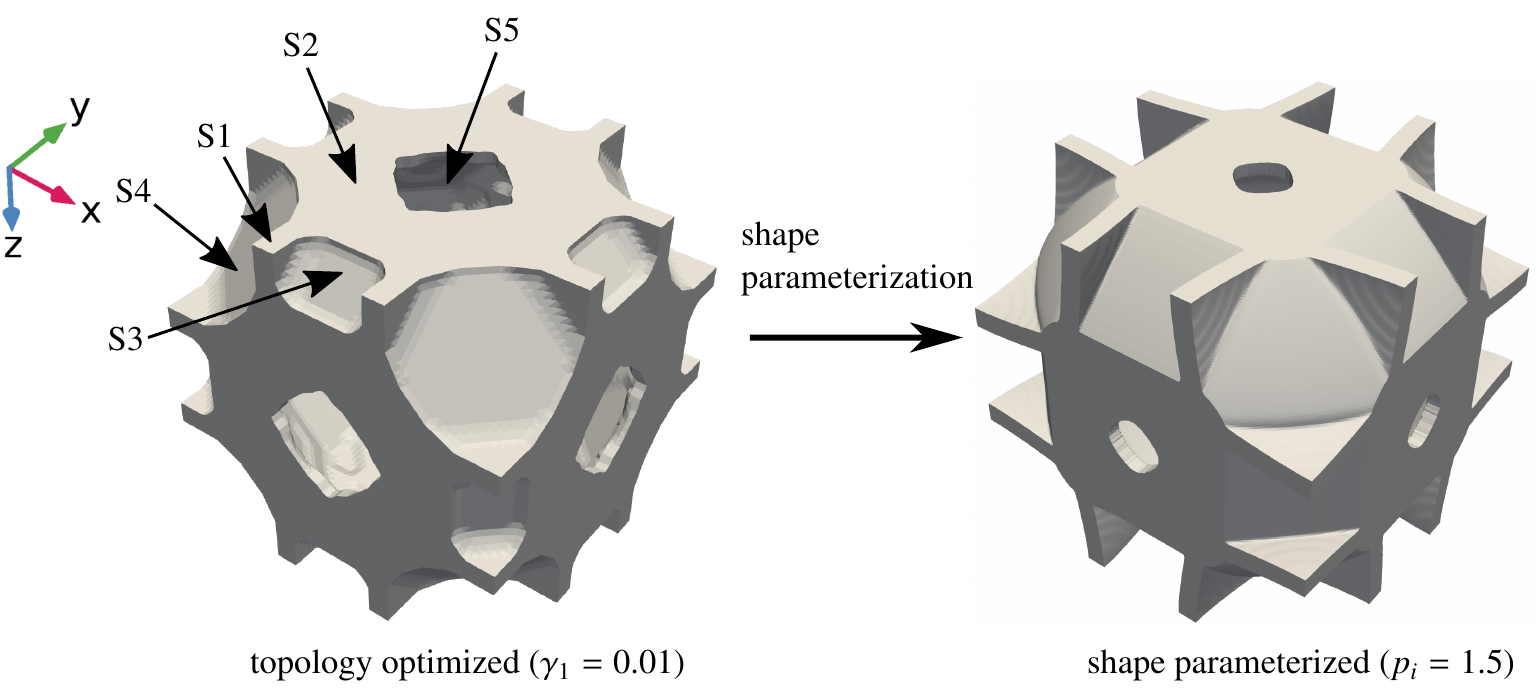} }\\ 
		\hline 
		Si & S1 &S2&S3 &S4 &  S5   \\
		\hline
		$\mathbf{c}_{i}$  &[0.0,0.0,0.0] &[0.5,0.5,0.5]  & [0.5,0.5,0.5] &[0.5,0.5,0.5] &[0.5,0.5,0]   \\
		$\mathbf{z}_{i}$  &$[x_1,y_1,z_1]$& $[0.5,0.5,z_2]$&$[0.5,y_3,y_3]$   &$[x_4,x_4,x_4]$ &$[x_5,x_5,0]$  \\
		$t_i$ & $t_1$&0.1 &  $t_3$ & $t_4$ &- \\ 
		$p_i$ & $p_1$& 1&$p_3$ &$p_4$ & $p_5$  \\ 
		\hline
	\end{tabular}	
	\caption{Illustration of shape parameterization based on the topology optimized microstructure  and  geometric design parameters of five control  super-ellipsoids. Top row:  Left: Topology optimized microstructure obtained with   $\gamma_1=0.01$, i.e., the topology optimized microstructure in the middle in Fig.~\ref{fig:OptDesUniShapeAllE}; Right: Shape parameterized microstructure with $p_i=1.5$ and $t_i=0.05$ and the detailed parameterization process is illustrated in Fig.~\ref{fig:ShapePara} as the case of $p_i=1.5$.  Symbols represent design variables  and `-' indicates  non-existing. The thickness for S2 is   preset to be a big value to ensure that it models a plate laying on the microstructure boundaries. } 	\label{Tab:Para}
\end{table}

The topology optimized microstructures denoted by the black line in Fig.~\ref{fig:OptDesUniShapeAllE}  can be parameterized using four hollow  (S1, S2, S3, S4) and one  solid   (S5) super-ellipsoids, as illustrated in the top row of  Tab.~\ref{Tab:Para}.  The geometric parameters associated with each one are listed in Tab.~\ref{Tab:Para}.  The detailed parameterization process is illustrated in Fig.~\ref{fig:ShapePara} (a).  In this study, we constrain the power to $p_i\leq 4$, and the location variables between 0.0 and 0.5. The control point of S1 can move freely in the 1/8 microstructure in order to increase the design freedom. The control point of S2  moves along the  z-direction and S2 is used to generate a plate normal to the z-direction.  The ones of S3 and S5  move along face diagonals. S3 is used to generate an extruded hollow super-ellipse along the x-direction and S5 is used to generate an extruded super-ellipse along the z-direction at the microstructure center.  The control point of  S4 moves on the body diagonal and S4 is a hollow super-ellipsoid with equal semi-diameters.   Hence, the microstructure architecture is controlled by 14  design parameters, i.e., $\left[ x_1 ,  y_1, z_1,  t_1, p_1, z_2, y_3, t_3, p_3, x_4,t_4, p_4, x_5, p_5 \right]$. 
 
Fig.~\ref{fig:ShapePara} illustrates how a 1/8 microstructure of   $[0,\ 0.5] \times [0,\ 0.5] \times [0,\ 0.5] $ is parameterized using the four hollow  and one solid super-ellipsoid with geometric cubic symmetry enforced.  As illustrated in Fig.~\ref{fig:ShapePara},  S2, S3, and S4 together generate a hollow geometry to represent the main body of microstructures. It has to be pointed out here that the inner and outer shapes of the hollow geometry are the intersection  of the inner  and outer super-ellipsoids of S2, S3  and S4, respectively.  S1 represents the  members connecting to the microstructure's main body which can change from flat to curved plates as shown in Fig.~\ref{fig:ShapePara}.   S5 is used to punch holes at each face center.  Fig.~\ref{fig:ShapePara}~(a), (b) and (c) display the generation process of microstructures with $p_i=1.5$,  $p_i=1 $  and  $p_i=2$, respectively.
 
The final design variables in the microstructure are numerically calculated using the density variables generated using the five features, written as
\begin{equation}\label{eq:SupFin}
{\rho}_e =  \left( 1- \left( 1-{\rho}^{2,1}_e{\rho}^{3,1}_e{\rho}^{4,1}_e \right) \left( 1-{\rho}^{1,1}_e\left( 1-{\rho}^{10}_e\right) \right) \right) \left( 1-{\rho}^{2,0}_e{\rho}^{3,0}_e{\rho}^{4,0}_e\right)\left( 1-{\rho}^{5,0}_e\right)  
\end{equation}

To smooth the intersection between different features, the density filter is employed with $r=0.025$ together with the projection in Eq. \eqref{Eq:2DProj} with $\beta_1=50$. The sensitivities of the physical density $\bar{\rho}_e$ with respect to each design variable  can be analytically derived using Eq. \eqref{eq:TPre}-\eqref{eq:SupFin}. The sensitivities of objective and constraints with respect to design variables are obtained using the chain rule,  combining the sensitivities of objective and constraints with respect to the physical densities and the sensitivities of the physical densities to the design variables. MMA is employed to update the designs.

\begin{figure}[!htb]
	\centering
 \includegraphics[]{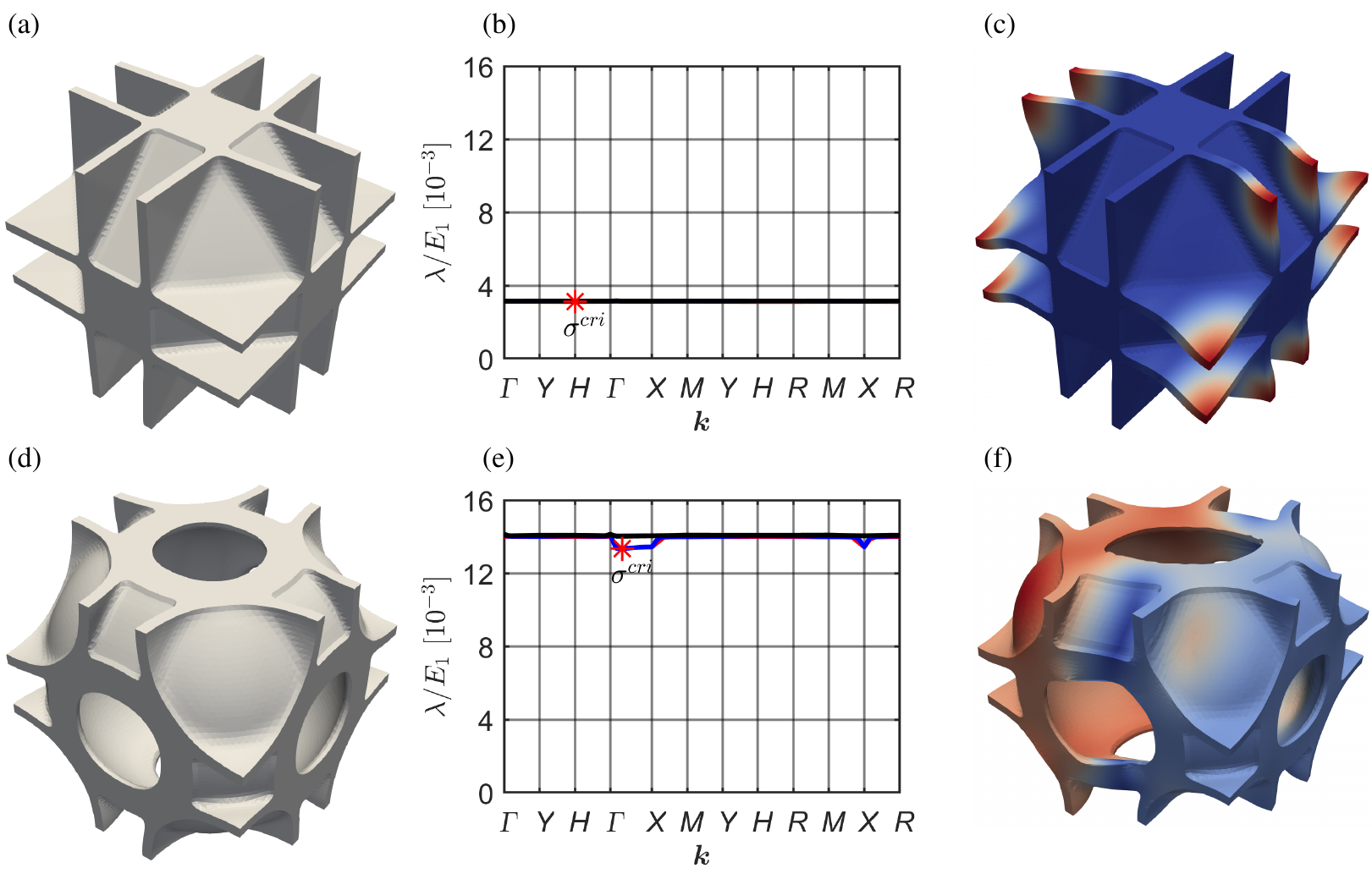}
	\caption{Shape optimized microstructures. (a)-(c): Shape optimized  microstructure  for maximum  $\bar{E}$,  ${\bar{E}}/E_1=0.0864$. (d)-(f) Shape optimized microstructure for maximum uniaxial buckling strength, ${\bar{E}}/E_1=0.0640$ .  Left column: Unit cell; Middle column: band diagram; Right column: Critical buckling mode.   }
	\label{fig:OptDesShape}
\end{figure}

Fig.~\ref{fig:OptDesShape} presents two shape optimized isotropic microstructures and corresponding band diagrams and critical buckling modes, i.e., microstructures for maximum  $\bar{E}$, see   (a)-(c), and   maximum $\sigma^{cri}$  see (d)-(f).  The shape optimized microstructure for maximum $\bar{E}$ buckles at the H point, and its first three buckling bands are tightly stacked.  The Young's-modulus-optimized microstructure exhibits the lowest buckling strength with  ${\sigma^{cri}}/{E_1}=0.00313$, which is lower than the hollow sphere in Fig.~\ref{fig:Illustration}. But the optimized microstructure exhibits an improvement of $173\%$ in Young's modulus and $272\%$ in the bulk modulus compared to the  hollow sphere.   The shape optimized microstructure for maximum ${\sigma^{cri}}$ buckles in a global shear mode. Due to the geometrical restriction imposed by the super-ellipsoids, bigger variations are observed in the first buckling band compared to the corresponding topology optimized microstructure in Fig.~\ref{fig:OptDes} (d)-(f). Hence the shape-optimized material exhibits decreased buckling strength with ${\sigma^{cri}}/E_1=0.01335$ and stiffness  with ${\bar{E}}/E_1=0.0640$  and ${\bar{\kappa}}/E_1=0.0504$. 
 
A set of shape optimized microstructures with tunable stiffness and buckling strength are again obtained by choosing  $\gamma_1$ to 0, 0.001, 0.01, 0.1, and 1. The corresponding performance is presented in Fig.~\ref{fig:OptDesUniShapeAllE}, represented by the red line with circles indicating each shape optimized microstructure. Compared to the topology optimized  microstructures, geometrical simplifications introduced by the shape optimization result in slightly decreased stiffness and buckling strength.  The maximum buckling strength obtained by the shape optimized  microstructures is $89\%$ of the corresponding topology optimized one, while the maximum Young's modulus is  $98\%$ of the corresponding topology optimized microstructure for the maximum Young's modulus. 

Compared to the isotropic stiffness-optimal plate microstructure in Fig.~\ref{fig:OptDesUniShapeAllE},  the topology optimized microstructure class covers Young's modulus range from   $81\%$    to $62\%$, a bulk modulus range from $100\%$ to $78\%$, and a buckling strength range from $422\%$  to  $ 863\%$. The simplified microstructures obtained from shape optimization cover a Young's modulus range from   $79\%$    to $58\%$, a bulk modulus range from $100\%$ to $67\%$, and a buckling strength range from  $180\%$  to  $767\%$.  

All the optimized microstructures have smooth features, and there is no extreme stress concentration in the designs. Compared to the isotropic stiffness-optimal plate material, the optimized microstructures possess smaller stiffness, hence, undergo larger displacement under the same macroscopic stress situation, which leads to smaller yield strength. The yield strength of the optimized microstructures monotonically decreases as the microstructure stiffness decreases.	

A recent study, Ref.~\cite{Andersen2021}, shows that stiffness-isotropic materials typically do not show buckling strength isotropy.  The buckling strength of the shape-optimized microstructure  with maximum buckling strength in Fig.~\ref{fig:OptDesShape} (d)  is evaluated for uniaxial compression along with two other directions,  i.e., body and face diagonals. The corresponding relative buckling strengths are 0.01364  and 0.01047. Hence, the optimized microstructures in this study do not possess buckling strength isotropy. Furthermore, it is shown in Ref~\cite{Andersen2020} that, in terms of buckling strength,  the volume fraction transition between simple isotropic stiffness-optimal truss and plate microstructures is around 31\%.   However, our optimization results show that this transition volume fraction can be lowered by the obtained hybrids between plate and truss  microstructures. Hence, for high volume fractions the pure plate structures will be best. For very low volume fractions (and an imposed minimum length scale), only viable structures are simple truss microstructures. Our optimizations thus point to novel microstructures that provide the best transition between the two extreme states.

 \section{Conclusion}
 3D isotropic microstructures with tunable stiffness and buckling strength have  systematically been designed. The effective material properties are evaluated using homogenization, and linear buckling analysis is employed to predict the material buckling strength for a given macroscopic stress state, where both microscopic and macroscopic failure modes are captured using  Bloch-Floquet boundary conditions. The optimization problem is formulated to improve the weighted stiffness and bucking strength via redistribution of base material in the material microstructure.   The one-case robust topology optimization formulation is employed to ensure single length-scale material design. 

The proposed optimization formulation is validated by designing a microstructure for maximum buckling strength under hydrostatic compression. Numerical results show that length scale is imposed in both solid and void regions and that the imposed length scale restrains the evolution of material geometrical hierarchy for the given volume fraction of $20\%$, and favors single length-scale material microstructures. 3D isotropic microstructures with tunable stiffness and buckling strength under uniaxial compression are obtained using the proposed formulation. Detailed comparisons between  designs obtained with different weighting between stiffness and buckling strength demonstrate  the evolution  from flat plates for high stiffness  to plate-truss hybrids with  increasing member thickness can increase the bending stiffness of microstructures  and hence  significantly enhances  material buckling strength at a small stiffness cost.

Inspired by the topology optimized material architectures, a subsequent shape optimization is proposed by parametrizing material configurations using a set of hollow and a solid super-ellipsoids. The material configurations are further optimized by using several design parameters hence achieving tunable stiffness and buckling strength. Due to geometrical restrictions introduced by the shape optimization, the shape optimized materials exhibit somewhat inferior stiffness and buckling strength compared to the corresponding topology optimized materials, but on the other hand are geometrically simpler. Compared to the isotropic stiffness-optimal plate material from~\cite{Berger2017}, the simplified material class covers a Young's modulus range from   $79\%$    to $58\%$, a bulk modulus range from $100\%$ to $67\%$, and a buckling strength range from  $180\%$  to  $767\%$ of that reference microstucture.  Moreover, buckling strength of the optimized microstructures shows a weak dependence on the wave vector. 

The optimized microstructures represent the stiffest and strongest architected materials to date. Experimental verifications are on-going. Even though this study focuses on enhancing buckling strength under uniaxial and hydrostatic compression, the proposed optimization formulation applies to any given macroscopic stress situation.

\section*{Acknowledgments}
We acknowledge the financial support from the Villum Fonden  through the Villum Investigator Project InnoTop. We further acknowledge valuable discussions with Morten N. Andersen at the Department of Mechanical Engineering in Technical University of Denmark.


\end{document}